\documentclass[12pt]{article}
\usepackage{latexsym}
\usepackage{amssymb}
\usepackage{graphicx}
\usepackage{cite}

\newtheorem{Theorem}{Theorem}[part]
\newtheorem{Definition}{Definition}[part]
\newtheorem{Proposition}{Proposition}[part]

\newtheorem{Corollary}{Corollary}[part]
\newtheorem{Remark}{Remark}[part]
\newtheorem{Example}{Example}[part]

\newtheorem{Condition}{Condition}[part]

\def \ep{\hbox{ }\hfill$\Box$}

\begin{document}
\title{Exceptionally Regular Tensors and Tensor Complementarity Problems
\thanks{This work was partially supported by the National Natural
Science Foundation of China (Grant No. 11431002).}
}

\author{Yong Wang\thanks{Department of Mathematics, School of Science, Tianjin University, Tianjin 300072, P.R. China. Email: wang\_yong@tju.edu.cn. }
\and
Zheng-Hai Huang\thanks{Corresponding Author. Department of Mathematics, School of Science, Tianjin University, Tianjin 300072, P.R. China. This author is also with the Center for Applied Mathematics of Tianjin University. Email: huangzhenghai@tju.edu.cn. Tel:+86-22-27403615 Fax:+86-22-27403615}
\and
Xue-Li Bai\thanks{Department of Mathematics, School of Science, Tianjin University, Tianjin 300072, P.R. China. Email: bai7707@163.com}
}

\date{August 26, 2015}

\maketitle

\begin{abstract}
\noindent
Recently, many structured tensors are defined and their properties are discussed in the literature. In this paper, we introduce a new class of structured tensors, called exceptionally regular tensor, which is relevant to the tensor complementarity problem. We show that this class of tensors is a wide class of tensors which includes many important structured tensors as its special cases. By constructing two examples, we demonstrate that an exceptionally regular tensor can be, but not always, an $R$-tensor. We also show that within the class of the semi-positive tensors, the class of exceptionally regular tensors coincides with the class of $R$-tensors. In addition, we consider the tensor complementarity problem with an exceptionally regular tensor or an $R$-tensor or a $P_0+R_0$-tensor, and show that the solution sets of these classes of tensor complementarity problems are nonempty and compact. \vspace{3mm}

\noindent {\bf Key words:}\hspace{2mm} Tensor complementarity problem, exceptionally regular tensor, weak $P$-tensor, $R$-tensor, $Z$-eigenvalue, $H$-eigenvalue \vspace{3mm}

\noindent {\bf Mathematics Subject Classifications(2000):}\hspace{2mm} 90C33, 65K10, 15A18, 15A69, 65F15, 65F10. \vspace{3mm}

\end{abstract}
\section{Introduction}
\setcounter{equation}{0} \setcounter{Assumption}{0}
\setcounter{Theorem}{0} \setcounter{Proposition}{0}
\setcounter{Corollary}{0} \setcounter{Lemma}{0}
\setcounter{Definition}{0} \setcounter{Remark}{0}
\setcounter{Algorithm}{0}

\hspace{4mm} A real $m$-order $n$-dimensional tensor ${\cal A}$ with $m,n$ being positive integers is an $m$-way array, which can be denoted by ${\cal A}=(a_{i_1i_2\cdots i_m})$ with $a_{i_1i_2\cdots i_m}\in \mathbb{R}$ for all $i_j\in \{1,2,\ldots,n\}$ and $j\in \{1,2,\ldots,m\}$. Obviously, ${\cal A}$ is a matrix if $m=2$. Throughout this paper, we assume that $m\geq 3$ and $n\geq 2$ unless otherwise stated; and use $\mathbb{T}_{m,n}$ to denote the set of all real $m$-order $n$-dimensional tensors. As a generalization of matrix, tensor has attracted more and more attention in the last ten years.

As the generalization of eigenvalues of matrices, eigenvalues of tensors are initially introduced and studied by Qi \cite{q-05} and Lim \cite{l-05}. They have been applied in many aspects, including magnetic resonance imaging \cite{hhnq-12, qyw-10} and hypergraph theory \cite{hq-12, lqy-13}. Several different concepts of eigenvalues of tensors are given in the literature. Two popular concepts are the $H$-eigenvalue and the $Z$-eigenvalue. Recall that for any given ${\cal A}\in \mathbb{T}_{m,n}$ and $x\in\mathbb{R}^n$, ${\cal A}x^{m-1}$ denotes a vector in $\mathbb{R}^n$ and
$$
\left({\cal A}x^{m-1}\right)_i:=\sum_{i_2,\cdots,i_m=1}^na_{ii_2\cdots i_m}x_{i_2}\ldots x_{i_m},\; \quad \forall i\in \{1,2,\ldots,n\}.
$$
If there exist a nonzero vector $x\in \mathbb{R}^n$ and a scalar $\lambda\in \mathbb{R}$ such that
$$
\left({\cal A}x^{m-1}\right)_i=\lambda x_i^{m-1},\;\forall i\in \{1,2,\ldots,n\},
$$
then $\lambda$ is called an $H$-eigenvalue of ${\cal A}$ and $x$ is called an $H$-eigenvector of ${\cal A}$ associated with $\lambda$. If there exist a nonzero vector $x\in \mathbb{R}^n$ and a scalar $\lambda\in \mathbb{R}$ such that
$$
\left({\cal A}x^{m-1}\right)_i=\lambda x_i,\;\forall i\in \{1,2,\ldots,n\} \quad\mbox{\rm and}\quad x^Tx=1,
$$
then $\lambda$ is called a $Z$-eigenvalue of ${\cal A}$ and $x$ is called a $Z$-eigenvector of ${\cal A}$ associated with $\lambda$.

Recall that for any given $A\in \mathbb{R}^{n\times n}$ and $q\in \mathbb{R}^n$, the linear complementarity problem, denoted by LCP$(q,A)$, is to find a point $x \in \mathbb{R}^n$ such that
\begin{eqnarray*}
x\geq 0,\quad Ax+q\geq 0,\quad x^T(Ax+q)=0.
\end{eqnarray*}
As a generalization of LCP$(q,A)$, the tensor complementarity problem has been introduced and investigated \cite{sq-14,sq-15,lqx-15}. For any given ${\cal A}\in \mathbb{T}_{m,n}$ and $q\in \mathbb{R}^n$, the tensor complementarity problem, denoted by TCP$(q,{\cal A})$, is to find a point $x \in \mathbb{R}^n$ such that
\begin{eqnarray*}
x\geq 0,\quad {\cal A}x^{m-1}+q\geq 0,\quad x^T({\cal A}x^{m-1}+q)=0.
\end{eqnarray*}
It is well known that properties of various structured matrices play important roles in theory and methods for LCP$(q,A)$. There is no doubt that properties of various structured tensors will play important roles in the study of TCP$(q,{\cal A})$. Recently, many structured tensors are introduced, such as $M$-tensor \cite{zqz-14,dqw-13}, $P\;(P_0$)-tensor \cite{sq-14}, $B\;(B_0$)-tensor \cite{sq-14,qs-14}, and $R\;(R_0$)-tensor;  and their properties are studied.

Recall that for any given nonlinear function $f: \mathbb{R}^n\rightarrow \mathbb{R}^n$, the nonlinear complementarity problem, denoted by CP$(f)$, is to find a point $x \in \mathbb{R}^n$ such that
\begin{eqnarray*}
x\geq 0,\quad f(x)\geq 0,\quad x^Tf(x)=0.
\end{eqnarray*}
It is obvious that TCP$(q,{\cal A})$ is a special class of CP$(f)$ with $f(x)={\cal A}x^{m-1}+q$.
It is well known that the concept of {\it exceptional family of elements} is a powerful tool for investigating CP$(f)$, and many good theoretical results for CP$(f)$ are obtained by using such a tool \cite{isac-97,isac-00,isac-98,zi-00,zh-99,zhq-99,hhxz-01}. In \cite{zi-00}, Zhao and Isac introduced the notion of {\it exceptionally regular function}, and  based on this notion, they established an existence theorem for CP$(f)$ with $f$ being a positively homogeneous function.

Inspired by these works, we will introduce a new class of structured tensors with the help of exceptionally regular function, which is called {\it exceptionally regular tensor}; and study related properties of this class of tensors.

The rest of this paper is organized as follows. In the next section, we review several important classes of structured tensors and some known related properties. We also discuss some properties of $Z$-eigenvalues of the weak $P$-tensor, recently defined by Ding, Luo, and Qi \cite{dlq-15}. In Section 3, we define a new class of structured tensors, exceptionally regular tensor, and discuss some of its properties. We show that a weak $P$-tensor must be an exceptionally regular tensor, but the converse does not hold. We also show that $ER$-tensors and $R$-tensors are two different classes of tensors although their intersection is nonempty. In Section 4, we consider the tensor complementarity problem with an exceptionally regular tensor, and show that its solution set is nonempty and compact, which is an extension of the results obtained by \cite[Theorem 4.5]{cqw} and \cite[Theorem 6.2]{dlq-15}. We also show that the solution set of the tensor complementarity problem, with an $R$-tensor or a $P_0+R_0$-tensor, is nonempty and compact.  The final conclusions are given in Section 5.

Throughout this paper, for any positive integer $n$, we denote $[n]:=\{1,2,\ldots,n\}$ and $\mathbb{R}^n_+:=\{x\in \mathbb{R}^n: x\geq 0\}$.

\section{Preliminaries}
\setcounter{equation}{0} \setcounter{Assumption}{0}
\setcounter{Theorem}{0} \setcounter{Proposition}{0}
\setcounter{Corollary}{0} \setcounter{Lemma}{0}
\setcounter{Definition}{0} \setcounter{Remark}{0}
\setcounter{Algorithm}{0}

\hspace{4mm}
In this section, we review definitions and properties of several structured tensors, which are useful for our sequential discussions. We also discuss the properties of eigenvalues of two classes of tensors.

\begin{Definition} \label{rdef} A tensor $\mathcal{A}\in \mathbb{T}_{m,n}$ is said to be
\begin{itemize}
  \item[(i)] {\bf semi-positive}\cite{sq-15} iff for each $x\in \mathbb{R}^n_+\setminus\{0\}$, there exists an index $i\in [n]$ such that
  $$x_i>0 \quad \mbox{\rm and}\quad (\mathcal{A}x^{m-1})_i\geq 0;$$
  \item[(ii)] {\bf strictly semi-positive}\cite{sq-15} iff for each $x\in \mathbb{R}^n_+\setminus\{0\}$, there exists an index $i\in [n]$ such that
  $$x_i>0 \quad \mbox{\rm and}\quad (\mathcal{A}x^{m-1})_i> 0;$$
  \item[(iii)] {\bf a $P_0$-tensor}\cite{sq-14} iff for each $x\in \mathbb{R}^n\setminus\{0\}$, there exists an index $i\in [n]$ such that
  $$x_i\neq 0\quad \mbox{\rm and}\quad x_i(\mathcal{A}x^{m-1})_i\geq0;$$
  \item[(iv)] {\bf a $P$-tensor}\cite{sq-14} iff for each $x\in \mathbb{R}^n\setminus\{0\}$, there exists an index $i\in [n]$ such that $$x_i(\mathcal{A}x^{m-1})_i>0;$$
  \item[(vi)] {\bf strictly copositive}\cite{sq-150} iff $\mathcal{A}x^m> 0$ for all $x\in \mathbb{R}_+^n\setminus\{0\}$;
  \item[(vii)] {\bf positive define}\cite{sq-150} iff $\mathcal{A}x^m>0$ for all $x\in \mathbb{R}^n\setminus\{0\}$;
  \item[(viii)]{\bf a $Q$-tensor} \cite{sq-15} iff for any $q\in \mathbb{R}^n$, TCP$(q, \mathcal{A})$ has a solution.
\end{itemize}
\end{Definition}

Clearly, every strictly semi-positive tensor is a semi-positive tensor, and every $P_0$-tensor is certainly semi-positive.

About the $Z$-eigenvalue and the $H$-eigenvalue of a (strictly) semi-positive tensor, we have the following observation.
\begin{Theorem}
If $\mathcal{A}\in \mathbb{T}_{m,n}$ is (strictly) semi-positive, then the $Z$-eigenvalue associated with a nonnegative $Z$-eigenvector of $\mathcal{A}$ is nonnegative (positive); and the $H$-eigenvalue associated with a nonnegative $H$-eigenvector of $\mathcal{A}$ is nonnegative (positive).
\end{Theorem}
{\bf Proof.} We first show the result for the $Z$-eigenvalue holds. Suppose that $\lambda$ is a $Z$-eigenvalue of $\mathcal{A}$ associated with a nonnegative $Z$-eigenvector $\hat{x}$, i.e., $\mathcal{A} \hat{x}^{m-1}=\lambda \hat{x}$. Since $\mathcal{A}$ is (strictly) semi-positive and $\hat{x}\geq 0$ with $\hat{x}\neq 0$, there exists an index $i_0\in [n]$ such that $\hat{x}_{i_0}>0$ and $\lambda \hat{x}_{i_0}=(\mathcal{A} \hat{x}^{m-1})_{i_0}\geq 0(>0)$. So, we get $\lambda\geq0(>0)$.

Similarly, we can show that the $H$-eigenvalue associated with a nonnegative $H$-eigenvector of a (strictly) semi-positive tensor $\mathcal{A}$ is nonnegative (positive).\ep

\begin{Definition}\label{r0def} \cite{sq-15} A tensor $\mathcal{A}\in \mathbb{T}_{m,n}$ is called an $R$-tensor, if there exists no $(x, t)\in(\mathbb{R}^n_+\setminus \{0\})\times\mathbb{R}_+$ such that
\begin{equation}\label{rtensor}
\left\{\begin{array}{ll}
(\mathcal{A}x^{m-1})_i+t=0,\quad & \mbox{\rm if}\; x_i>0,\\
(\mathcal{A}x^{m-1})_i+t\geq0,\quad & \mbox{\rm if}\; x_i=0.
\end{array}\right.
\end{equation}
A tensor $\mathcal{A}\in \mathbb{T}_{m,n}$ is called an $R_0$-tensor, if the system (\ref{rtensor}) has no nonzero solution when $t=0$, i.e., there exists no $x\in \mathbb{R}^n_+\setminus \{0\}$ such that
\begin{equation}\label{r0tensor}
\left\{\begin{array}{ll}
(\mathcal{A}x^{m-1})_i=0,\quad & \mbox{\rm if}\; x_i>0,\\
(\mathcal{A}x^{m-1})_i\geq0,\quad & \mbox{\rm if}\; x_i=0.
\end{array}\right.
\end{equation}
\end{Definition}

It is obvious that every $R$-tensor is an $R_0$-tensor, but the converse does not hold. An example is given in \cite{sq-15}. Moreover, it is also demonstrated that if an $R_0$-tensor is also semi-positive, then this tensor is an $R$-tensor \cite[Theorem 3.4]{sq-15}.

 The following definition and theorem can be found in \cite{isac-97}.
\begin{Definition}\label{def-ex}
A set of points $\{x^k\}\subset \mathbb{R}^n_+$ is an exceptional family of elements for the continuous function $f$ if $\|x^k\|\rightarrow \infty$ as $k\rightarrow \infty$ and, for each $k>0$, there exists a scalar $\mu_k>0$ such that
\begin{equation}\label{E-ex-f}
\left\{\begin{array}{ll}
f_i(x^k)=-\mu_kx_i^k,\quad & \mbox{\rm if}\; x_i^k>0,\\
f_i(x^k)\geq0,\quad & \mbox{\rm if}\; x_i^k=0.
\end{array}\right.
\end{equation}
\end{Definition}

\begin{Theorem}\label{thm-ex-basic}
For any continuous function $f: \mathbb{R}^{n}_+\rightarrow \mathbb{R}^{n}$, there exists either a solution to CP$(f)$ or an exceptional family of elements for $f$.
\end{Theorem}

More recently, Ding, Luo and Qi\cite{dlq-15} defined a new class of $P$-type tensors:
\begin{Definition}\label{wpdef}
$\mathcal{A}\in \mathbb{T}_{m,n}$ is called a weak $P$-tensor, if for each nonzero $x\in\mathbb{R}^n$, there exists an index $i\in [n]$ such that
$x_i^{m-1} (\mathcal{A}x^{m-1})_i > 0.$
\end{Definition}

It is easy to see that when $m$ is even, the weak $P$-tensor defined by Definition \ref{wpdef} is consistent with the $P$-tensor defined by Definition \ref{rdef}. It has been known that when $m$ is odd, the $P$-tensor does not exist \cite{yy-14}. So, the weak $P$-tensor is an extension of the $P$-tensor. In \cite{dlq-15}, Ding, Luo and Qi also called the weak $P$-tensor as the $P$-tensor. In order to distinguish these two classes of different $P$-tensors, we call the $P$-tensor defined by Ding, Luo and Qi\cite{dlq-15} as the weak $P$-tensor ($wP$-tensor for short) in this paper.

In \cite{dlq-15}, many nice properties of $wP$-tensors were studied, including the properties of the $H$-eigenvalue. In the following, we consider the $Z$-eigenvalue of a $wP$-tensor.
\begin{Theorem}
Let $\mathcal{A}\in \mathbb{T}_{m,n}$ be a $wP$-tensor. When $m$ is even, all of its $Z$-eigenvalues are all positive; and when $m$ is odd, every of its $Z$-eigenvalues associated with a nonnegative (nonpositive) $Z$-eigenvector is positive (negative).
\end{Theorem}
{\bf Proof.} When $m$ is even, since the $wP$-tensor is just the $P$-tensor, the first result holds from \cite[Theorem 3.2]{sq-14}. In the following, we assume that $m$ is odd and show the second result. Suppose that $\lambda$ and $\hat{x}$ are a $Z$-eigenvalue and a corresponding $Z$-eigenvector of $\mathcal{A}$, respectively, i.e., $\mathcal{A} \hat{x}^{m-1}=\lambda \hat{x}$. Since $\mathcal{A}$ is a $wP$-tensor and $\hat{x}\neq 0$, there exists an index $i_0\in [n]$ such that $$\hat{x}_{i_0}^{m-1}(\mathcal{A} \hat{x}^{m-1})_{i_0}>0.$$
So, it follows that $\hat{x}_{i_0}\neq 0$.
\begin{itemize}
  \item If $\hat{x}\geq 0$, then $\hat{x}_{i_0}^{m-1}>0$, and hence, $(\mathcal{A} \hat{x}^{m-1})_{i_0}>0$. Furthermore, we get that $\lambda \hat{x}_{i_0}=(\mathcal{A} \hat{x}^{m-1})_{i_0}>0$. So, $\lambda>0$.
  \item If $\hat{x}\leq 0$, then $\hat{x}_{i_0}<0$, and hence, $\hat{x}_{i_0}^{m-1}>0$ since $m$ is odd. Thus, $(\mathcal{A} \hat{x}^{m-1})_{i_0}>0$. Furthermore, we get that $\lambda \hat{x}_{i_0}>0$. So, $\lambda<0$.
\end{itemize}
The proof is complete.\ep

\section{Exceptionally regular tensor}
\setcounter{equation}{0} \setcounter{Assumption}{0}
\setcounter{Theorem}{0} \setcounter{Proposition}{0}
\setcounter{Corollary}{0} \setcounter{Lemma}{0}
\setcounter{Definition}{0} \setcounter{Remark}{0}
\setcounter{Algorithm}{0}

\hspace{4mm}
In this section, we introduce the exceptionally regular tensor and discuss its properties.

The following concept can be found in \cite{zi-00}.
\begin{Definition} \cite{zi-00} The function $g(x)=f(x)-f(0)$ is exceptionally regular if there exists no $(x,\alpha)\in \mathbb{R}^{n}_+\times\mathbb{R}_+$ with $\|x\|_2=1$ such that
\begin{equation}\label{exreg}
\left\{\begin{array}{lll}
g_i(x)/x_i&=-\alpha,\quad & \mbox{\rm if}\; x_i>0,\\
g_i(x)    &\geq0,\quad & \mbox{\rm if}\; x_i=0.
\end{array}\right.
\end{equation}
\end{Definition}

Motivated by the concept of exceptionally regular function, we define a new class of structured tensors, which is called exceptionally regular tensor ($ER$-tensor for short).
\begin{Definition}\label{erdef}
$\mathcal{A}\in \mathbb{T}_{m,n}$ is called an $ER$-tensor, if there exists no $(x, t)\in(\mathbb{R}^n_+\setminus \{0\})\times\mathbb{R}_+$ such that
\begin{equation}\label{ertensor}
\left\{\begin{array}{lcll}
(\mathcal{A}x^{m-1})_i+tx_i&=&0,   \quad & \mbox{\rm if}\; x_i>0,\\
(\mathcal{A}x^{m-1})_i     &\geq&0,\quad & \mbox{\rm if}\; x_i=0.
\end{array}\right.
\end{equation}
\end{Definition}

When $m=2$, an $ER$-tensor reduces to a matrix, and we call it {\it an $ER$-matrix}. Thus, by Definition \ref{erdef}, we actually introduce a class of structured matrices, which is new to the best of our knowledge.

We first discuss the relationship between $ER$-tensor and related tensors.

From the definition of strictly semi-positive tensor, we can get that the class of strictly semi-positive tensors is a subset of the class of $ER$-tensors.
\begin{Proposition}\label{pro2}
If $\mathcal{A}\in \mathbb{T}_{m,n}$ is a strictly semi-positive tensor, then $\mathcal{A}$ is an $ER$-tensor.
\end{Proposition}
{\bf Proof.} Since $\mathcal{A}\in \mathbb{T}_{m,n}$ is strictly semi-positive, it follows that for any $x\in\mathbb{R}^n_+\setminus\{0\}$, there exists an index $i_0\in[n]$ such that $x_{i_0}>0$ and $(\mathcal{A}x^{m-1})_{i_0}>0$. That is to say, for any $x\in\mathbb{R}^n_+\setminus\{0\}$, the system (\ref{ertensor}) has no solution. So, $\mathcal{A}$ is an $ER$-tensor.\ep

\begin{Remark}
If $\mathcal{A}\in \mathbb{T}_{m,n}$ is a positive definite tensor, then for any $x\in\mathbb{R}^n\setminus\{0\}$, $x^T(\mathcal{A}x^{m-1})>0$. So, there exists at least one index $i\in [n]$ such that $x_i(\mathcal{A}x^{m-1})_i>0$. Therefore, $\mathcal{A}$ is a $P$-tensor. It is easy to see from Definition \ref{rdef}(ii)(iv) that a $P$-tensor is a strictly semi-positive tensor. And so is a positive tensor (all elements are positive). It is easy to see from Definition \ref{rdef}(i)(vi) that a strictly copositive tensor is also a strictly semi-positive tensor. Thus, it follows from Proposition \ref{pro2} that the tensors mentioned above are all $ER$-tensors.
\end{Remark}

\begin{Theorem}\label{wpiser}
If $\mathcal{A}$ is a $wP$-tensor defined by Definition \ref{wpdef}, then $\mathcal{A}$ is an $ER$-tensor.
\end{Theorem}
{\bf Proof.} Since $\mathcal{A} $ is a $wP$-tensor, it follows from Definition \ref{wpdef} that for any $x\in \mathbb{R}^n\setminus\{0\}$, there exists an index $i\in [n]$ such that $x_i^{m-1}(\mathcal{A} x^{m-1})_i>0$. In particular, for any $x\in\mathbb{R}^n_+\setminus\{0\}$, there exists an index $i_0$ such that
$$x_{i_0}^{m-1}(\mathcal{A} x^{m-1})_{i_0}>0.$$
So, $x_{i_0}\neq 0$, and then $x_{i_0}>0$. Furthermore, we get $(\mathcal{A} x^{m-1})_{i_0}>0$. That is to say, for any $x\in\mathbb{R}^n_+\setminus\{0\}$, there exists an index $i_0$ such that $x_{i_0}>0$ and $(\mathcal{A} x^{m-1})_{i_0}>0$, which implies that the system (\ref{ertensor}) has no solution. So, $\mathcal{A} $ is an $ER$-tensor.\ep
\begin{Remark}
Recall that the $wP$-tensor is defined by Definition \ref{wpdef}, and many classes of important structured tensors are the subclasses of $wP$-tensors, including positive definite tensors \cite[Proposition 3.1]{dlq-15}, strongly completely positive tensors \cite[Proposition 3.4]{dlq-15}, nonsingular $H$-tensors with all positive diagonal entries \cite[Proposition 4.1]{dlq-15}, strictly diagonally dominant tensors with positive diagonal entries \cite[Corollary 4.2]{dlq-15}, Cauchy tensors with mutually distinct entries of generating vector\cite[Corollary 4.4]{dlq-15}, addition tensors of $wP$-tensors and completely positive tensors \cite[Theorem 4.5]{dlq-15}, odd-order $B$-tensors or symmetric even-order $B$-tensors \cite[Corollary 4.6]{dlq-15}, and so on. Thus, from Theorem \ref{wpiser}, the classes of tensors mentioned above are all $ER$-tensors.
\end{Remark}

Therefore, $ER$-tensors are a wide class of tensors which includes many important tensors as its special cases.

From Definitions \ref{r0def} and \ref{erdef}, it seems that the definitions of $ER$-tensor and $R$-tensor have some similarities. In the following, we discuss the relationship between these two classes of tensors.

We have known that a strictly semi-positive tensor is an $R$-tensor\cite{sq-15}. This and Proposition \ref{pro2} imply that the intersection of the class of $ER$-tensors and the class of $R$-tensors is nonempty. In the following, we construct two examples, which show that the class of $ER$-tensors is different from the class of $R$-tensors.
\begin{Example}\label{exa1}
Let $\mathcal{A}=(a_{ijk})\in \mathbb{T}_{3,2}$, where $a_{111}=-16, a_{122}=1, a_{211}=-17, a_{222}=1$ and all other elements of $\mathcal{A}$ are zeros, then $\mathcal{A}$ is an $R$-tensor, but not an $ER$-tensor.
\end{Example}

It is obvious that for any $x\in \mathbb{R}^2_+$,
\begin{eqnarray*}
\mathcal{A}x^2=\left(\begin{array}{c}-16x_1^2+x_2^2\\-17x_1^2+x_2^2\end{array}\right).
\end{eqnarray*}
In the following, we show that the results given in Example \ref{exa1} hold.

First, we show that $\mathcal{A}$ is an $R$-tensor. In fact,
\begin{itemize}
  \item if $x_1>0$, then $(\mathcal{A}x^2)_1+t=-16x_1^2+x_2^2+t=0$, i.e.,  $x_2^2=16x_1^2-t$, but $(\mathcal{A}x^2)_2+t=-17x_1^2+x_2^2+t=-x_1^2<0$;
  \item and if $x_2>0$, then $(\mathcal{A}x^2)_2+t=-17x_1^2+x_2^2+t=0$, i.e., $17x_1^2=x_2^2+t>0$, but
  $(\mathcal{A}x^2)_1+t=-16x_1^2+x_2^2+t=-16x_1^2+17x_1^2-t+t=x_1^2>0$.
\end{itemize}
Therefore, $\mathcal{A}$ is an $R$-tensor.

Second, we show that $\mathcal{A}$ is not an $ER$-tensor. We consider
\begin{eqnarray}
x_1>0, &\;\;& (\mathcal{A}x^2)_1+tx_1=-16x_1^2+x_2^2+tx_1=0; \label{E-exam1-1}\\
x_2>0, &\;\;& (\mathcal{A}x^2)_2+tx_2=-17x_1^2+x_2^2+tx_2=0. \label{E-exam1-2}
\end{eqnarray}
From (\ref{E-exam1-1}) it follows that $x_2^2=x_1(16x_1-t)$ and $16x_1\geq t$. These and (\ref{E-exam1-2}) imply that
$$
(\mathcal{A}x^2)_2+tx_2=-17x_1^2+x_2^2+tx_2=-x_1^2-tx_1+t\sqrt{16x_1^2-tx_1}=0.
$$
So, by using $x_1>0$, we further get
\begin{equation}\label{equ4}
x_1^3+2tx_1^2-15t^2x_1+t^3=0.
\end{equation}
Without lose of generality, we set $t=1$ and try to show that the equation
\begin{equation}\label{equ4-1}
f(z):=z^3+2z^2-15z+1=0
\end{equation}
has a root on the interval $(\frac{1}{16},+\infty)$.
It is obvious that $f(z)\rightarrow +\infty$ as $z\rightarrow +\infty$.
Let the derivation of $f(z)$ equal to zero, i.e.,
$$3z^2+4z-15=0,$$
we get $z=\frac{5}{3}>\frac{1}{16}$ and $f(\frac{5}{3})=-\frac{373}{27}<0$. So, the equation (\ref{equ4-1}) has a root $z^*\in(\frac{5}{3}, +\infty)\subset(\frac{1}{16}, +\infty)$. Thus, $(x_1,t):=(z^*,1)$ solves the equation (\ref{equ4}). Furthermore, take $x_2=\sqrt{16(z^*)^2-z^*}$, then $(\bar{x},\bar{t})\in (\mathbb{R}^2_+\setminus\{0\})\times \mathbb{R}_+$, with $\bar{x}=(z^*, \sqrt{16(z^*)^2-z^*})^T$ and $\bar{t}=1$, solves the system (\ref{ertensor}). Therefore, $\mathcal{A}$ is not an $ER$-tensor.\ep

\begin{Example}\label{exa2}
Let $\mathcal{A}=(a_{ijk})\in \mathbb{T}_{3,2}$, where $a_{111}=1, a_{122}=-1, a_{211}=2, a_{222}=-1$ and all other elements of $\mathcal{A}$ are zeros, then $\mathcal{A}$ is an $ER$-tensor, but not an $R$-tensor.
\end{Example}

It is obvious that for any $x\in \mathbb{R}^2_+$,
\begin{eqnarray*}
\mathcal{A}x^2=\left(\begin{array}{c}x_1^2-x_2^2\\2x_1^2-x_2^2\end{array}\right).
\end{eqnarray*}
In the following, we show that the results given in Example \ref{exa2} hold.

First, we show that $\mathcal{A}$ is an $ER$-tensor. We consider the following two cases.
\begin{itemize}
  \item [(C1)] If $x_1>0$, $(\mathcal{A}x^2)_1+tx_1=x_1^2-x_2^2+tx_1=0$. Then $$x_2^2=x_1^2+tx_1>0\quad  \mbox{\rm and}\quad  x_2-x_1=\frac{tx_1}{x_2+x_1}.$$
      Thus, $x_2>0$. We need to show that the equation
      \begin{equation}\label{equ3-0}
      (\mathcal{A}x^2)_2+tx_2=2x_1^2-x_2^2+tx_2=x_1^2+t(x_2-x_1)=x_1^2+t\frac{tx_1}{x_2+x_1}=0
      \end{equation}
      has no solution. In fact, suppose that the equation (\ref{equ3-0}) has a solution, then it follows from (\ref{equ3-0}) that $x_1^3+x_1^2x_2+t^2x_1=0$, which is impossible because $x_1>0$ and $x_2>0$.
  \item [(C2)] If $x_2>0$, $(\mathcal{A}x^2)_2+tx_2=2x_1^2-x_2^2+tx_2=0$. Then
  \begin{equation}\label{equ2}x_1^2=\frac{x_2^2-tx_2}{2}\geq 0.\end{equation}
  \begin{itemize}
    \item [a)]If $x_1=0$, then
  $$(\mathcal{A}x^2)_1+tx_1=x_1^2-x_2^2+tx_1=-x_2^2\leq 0,$$ which contradicts the condition that $x_2>0$.
      \item [b)]If $x_1>0$, then it follows from (\ref{equ2}) that $x_2>t$ and
 \begin{eqnarray*}
 (\mathcal{A}x^2)_1+tx_1=x_1^2-x_2^2+tx_1=\frac{x_2^2-tx_2}{2}-x_2^2+t\sqrt{\frac{x_2^2-tx_2}{2}}.
 \end{eqnarray*}
  Let $(\mathcal{A}x^2)_1+tx_1=0$, we derive a contradiction. It follows from $(\mathcal{A}x^2)_1+tx_1=0$ that
  $$\frac{x_2^2}{2}+\frac{tx_2}{2} =t\sqrt{\frac{x_2^2-tx_2}{2}},$$
  and hence,
  \begin{equation}\label{equ3}x_2^2(x_2+t)^2 =t^2(2x_2^2-2tx_2).\end{equation}
  Since $x_2>0$, (\ref{equ3}) can be simplified as
  $$x_2^3+2tx_2^2-t^2x_2+2t^3=0,$$
  which does not hold since $x_2^3>t^2x_2$ from $x_2>t$ and $t\geq 0$.
  \end{itemize}
\end{itemize}
Therefore, the system (\ref{ertensor}) has no solution $(x, t)\in(\mathbb{R}^n_+\setminus \{0\})\times\mathbb{R}_+$, which demonstrates that $\mathcal{A}$ is an $ER$-tensor.

Second, we show that $\mathcal{A}$ is not an $R$-tensor. In fact, take $\bar{x}_1=0, \bar{x}_2=a>0$ and $\bar{t}=a^2$. Then, it is easy to check that $(\bar{x},\bar{t})\in (\mathbb{R}^2_+\setminus\{0\})\times \mathbb{R}_+$ is a solution of (\ref{rtensor}). Therefore, $\mathcal{A}$ is not an $R$-tensor.\ep

In Theorem \ref{wpiser}, we have showed that every $wP$-tensor is an $ER$-tensor. In fact, the inverse does not hold.
\begin{Theorem}\label{wpiser-1}
The class of $wP$-tensors is a proper subset of the class of $ER$-tensors.
\end{Theorem}
{\bf Proof.} From Theorem \ref{wpiser}, we only need to show that there exists a tensor which is an $ER$-tensor, but not a $wP$-tensor. Suppose that $\mathcal{A}\in \mathbb{T}_{3,2}$ is given by Example \ref{exa2}, then $\mathcal{A}$ is an $ER$-tensor. Next, we show that $\mathcal{A}$ is not a $wP$-tensor. In fact, let $\bar{x}=(0,1)^T$ and we have that
\begin{eqnarray*}
\left\{\begin{array}{ll}
\bar{x}_1=0,\quad & \bar{x}_1^{2}(\mathcal{A}\bar{x}^{2})_1=0;\\
\bar{x}_2=1>0,\quad & \bar{x}_2^{2}(\mathcal{A}\bar{x}^{2})_2=1\cdot (0-1)=-1<0.
\end{array}\right.
\end{eqnarray*}
That is to say, for given $\bar{x}=(0,1)^T$, there exists no index $i\in \{1,2\}$ such that $\bar{x}_i^{2}(\mathcal{A}\bar{x}^{2})_i>0$. Thereby, the tensor $\mathcal{A}$ is not a $wP$-tensor.\ep

Now, we discuss the properties of $ER$-tensors.

The following proposition provides several necessary conditions for a tensor being an $ER$-tensor.
\begin{Proposition}\label{pro1}
Suppose that the tensor $\mathcal{A}\in \mathbb{T}_{m,n}$ is an $ER$-tensor. Then the following results hold.
\begin{itemize}
  \item [(i)] $\mathcal{A}$ is an $R_0$-tensor.
  \item [(ii)] Every principal sub-tensor of $\mathcal{A}$ is also an $ER$-tensor.
  \item [(iii)] The $Z$-eigenvalue of $\mathcal{A}$ associated with a nonnegative $Z$-eigenvector is positive.
\end{itemize}
\end{Proposition}
{\bf Proof.} (i) Since $\mathcal{A}$ is an $ER$-tensor, then any point $(x,t)\in (\mathbb{R}^n_+\setminus \{0\})\times \mathbb{R}_+$ is not a solution of the system (\ref{ertensor}). So, the system (\ref{ertensor}) has no  nonzero solution when $t=0$, that is, the system (\ref{r0tensor}) has no solution $x\in \mathbb{R}^n_+\setminus \{0\}$. Therefore, $\mathcal{A}$ is an $R_0$-tensor.

(ii) Let $J\subset [n]$ and $|J|=r\;(1\leq r\leq n)$, then $\mathcal{A}^J_{r}$ is one of the principal sub-tensors of $\mathcal{A}$. Suppose $\mathcal{A}^J_{r}$ is not an $ER$-tensor, then there exists a point $(x_J, t)\in (\mathbb{R}^{r}_+\setminus\{0\})\times\mathbb{R_+}$ satisfying the system (\ref{ertensor}). Define $\bar{x}\in \mathbb{R}^n_+\setminus\{0\}$ by
$$\bar{x}_i=\left\{\begin{array}{ll}
(x_J)_i,& \quad \mbox{\rm if}\; i\in J,\\
0,& \quad \mbox{\rm if}\; i\notin J,\\
\end{array}\right.$$
then, it is easy to see that the point $(\bar{x}, t)\in (\mathbb{R}^{n}_+\setminus\{0\})\times\mathbb{R_+}$ solves the system (\ref{ertensor}). Therefore, $\mathcal{A}$ is not an $ER$-tensor,  which causes a contradiction.

(iii) We assume that $x\geq 0$ is a $Z$-eigenvector of $\mathcal{A}$ and $\lambda$ is the corresponding $Z$-eigenvalue, then $\mathcal{A}x^{m-1}=\lambda x$. Suppose $\lambda\leq 0$ and let $t=-\lambda\geq 0$, then we have
\begin{eqnarray*}
\left\{\begin{array}{lcl}
(\mathcal{A}x^{m-1})_i+tx_i&=0, \quad & \mbox{\rm if}\;  x_i>0,\\
(\mathcal{A}x^{m-1})_i&\geq 0, \quad & \mbox{\rm if}\;  x_i=0.
\end{array}\right.
\end{eqnarray*}
That is, the point $(x, t)\in (\mathbb{R}^{n}_+\setminus\{0\})\times\mathbb{R_+}$ solves the system (\ref{ertensor}), which contradicts that $\mathcal{A}$ is an $ER$-tensor. Hence, we have $\lambda>0$. \ep

\begin{Corollary}\label{cor-1}
Given an $ER$-tensor ${\cal A}\in \mathbb{T}_{m,n}$ with $m$ being even, then the $Z$-eigenvalue of $\mathcal{A}$ associated with a nonnegative (or nonpositve) $Z$-eigenvector is positive.
\end{Corollary}

{\bf Proof.} By Proposition \ref{pro1}(iii), we only need to prove that the $Z$-eigenvalue of $\mathcal{A}$ associated with a nonpositve $Z$-eigenvector is positive.

Now we assume that $x\leq 0$ is a $Z$-eigenvector of $\mathcal{A}$ and $\lambda$ is the corresponding $Z$-eigenvalue, then $\mathcal{A}x^{m-1}=\lambda x$. Suppose $\lambda\leq 0$. Since $m$ is even, we have
$$
{\cal A}(-x)^{m-1}=-{\cal A}x^{m-1}=-\lambda x=\lambda (-x),
$$
which implies that $-x$ is also a $Z$-eigenvector associated with $\lambda$. Obviously, $-x\geq 0$. Let $t=-\lambda$ and $\bar{x}=-x$, then $(\bar{x}, t)\in (\mathbb{R}^{n}_+\setminus\{0\})\times\mathbb{R_+}$ solves the system (\ref{ertensor}), which contradicts that ${\cal A}$ is an $ER$-tensor. Therefore, we have $\lambda>0$. \ep

From Corollary \ref{cor-1} it follows that if an $ER$-matrix $A\in \mathbb{R}^{n\times n}$ has a nonpositive eigenvalue $\lambda$ associated with an eigenvector $x$, then there exist at least two distinct indexes $i,j\in [n]$ such that $x_ix_j<0$.

The following theorem gives the equivalence of three classes of structured tensors within the class of semi-positive tensors.
\begin{Theorem}\label{thm-equi}
If $\mathcal{A}\in \mathbb{T}_{m,n}$ is semi-positive, then the following results are equivalent.
\begin{itemize}
\item[(i)] $\mathcal{A}$ is an $R_0$-tensor,
\item[(ii)] $\mathcal{A}$ is an $ER$-tensor,
\item[(iii)] $\mathcal{A}$ is an $R$-tensor.
\end{itemize}
\end{Theorem}
{\bf Proof.} On one hand, it is obvious that every $R$-tensor is an $R_0$-tensor; and on the other hand, it follows from \cite[Theorem 3.4]{sq-15} that every semi-positive $R_0$-tensor is an $R$-tensor. Thus, (i) holds if and only if (iii) holds.

In the following, we show that (i) holds if and only if (ii) holds. Since every $ER$-tensor is an $R_0$-tensor by Proposition \ref{pro1}(i), we only need to show that $\mathcal{A}$ is an $ER$-tensor under the assumption that it is an $R_0$-tensor. Suppose that $\mathcal{A}$ is not an $ER$-tensor, then there exists a point $(\bar{x}, \bar{t})\in(\mathbb{R}^n_+\setminus \{0\})\times\mathbb{R}_+$ satisfying the system (\ref{ertensor}). Since $\mathcal{A}$ is an $R_0$-tensor, we have $\bar{t}>0$. Thus, we have
\begin{eqnarray*}
\left\{\begin{array}{lcll}
(\mathcal{A}\bar{x}^{m-1})_i+\bar{t}\bar{x}_i&=&0,   \quad & \mbox{\rm if}\; \bar{x}_i>0,\\
(\mathcal{A}\bar{x}^{m-1})_i     &\geq&0,\quad & \mbox{\rm if}\; \bar{x}_i=0,
\end{array}\right.
\end{eqnarray*}
i.e.,
\begin{eqnarray*}
\left\{\begin{array}{lcll}
(\mathcal{A}\bar{x}^{m-1})_i &=&-\bar{t}\bar{x}_i<0,   \quad & \mbox{\rm if}\; \bar{x}_i>0,\\
(\mathcal{A}\bar{x}^{m-1})_i &\geq&0,\quad & \mbox{\rm if}\; \bar{x}_i=0.
\end{array}\right.
\end{eqnarray*}
This implies that for $\bar{x}\in \mathbb{R}^n_+\setminus \{0\}$, we have
$$
\bar{x}_i\left({\cal A}\bar{x}^{m-1}\right)_i=-\bar{t}\bar{x}_i^2<0, \forall i\in \{j\in [n]: \bar{x}_j>0\},
$$
which contradicts the condition that $\mathcal{A}$ is a semi-positive tensor. So, $\mathcal{A}$ is an $ER$-tensor.

The proof is complete. \ep

Since every $P_0$-tensor is semi-positive, from Theorem \ref{thm-equi} we have the following results.
\begin{Corollary}\label{cor-thm-equi}
If $\mathcal{A}\in \mathbb{T}_{m,n}$ is a $P_0$-tensor, then $\mathcal{A}$ is an $R_0$-tensor iff $\mathcal{A}$ is an $ER$-tensor iff $\mathcal{A}$ is an $R$-tensor.
\end{Corollary}

\section{Properties of the solution set of TCP($q,\mathcal{A}$)}
\setcounter{equation}{0} \setcounter{Assumption}{0}
\setcounter{Theorem}{0} \setcounter{Proposition}{0}
\setcounter{Corollary}{0} \setcounter{Lemma}{0}
\setcounter{Definition}{0} \setcounter{Remark}{0}
\setcounter{Algorithm}{0}
\hspace{4mm}
In this section, we study properties of the solution set of TCP$(q, \mathcal{A})$. For any $x\in \mathbb{R}^n$, we denote
$$
[x]_+:=(\max\{x_1,0\},\ldots,\max\{x_n,0\})^T.
$$
We will use the following condition and proposition.
\begin{Condition}\label{cond-1}
Given $\mathcal{A}\in \mathbb{T}_{m,n}$ and $q\in \mathbb{R}^n$.
If there exists a sequence $\{x^k\}\subset\mathbb{R}^n_+$ satisfying
\begin{eqnarray}\label{E-cond-1}
\|x^k\|\rightarrow +\infty \quad \mbox{\rm and}\quad  \frac{[-\mathcal{A}(x^k)^{m-1}-q]_+}{\|x^k\|}\rightarrow 0 \;\mbox{as }\; k\rightarrow +\infty,
\end{eqnarray}
then there exists an index $i\in [n]$ such that $x_i[\mathcal{A}(x^k)^{m-1}+q]_i>0$ holds for some $k\geq 0$.
\end{Condition}

\begin{Proposition}\label{prop-basic-1}
Given $\mathcal{A}\in \mathbb{T}_{m,n}$ and $q\in \mathbb{R}^n$.
\begin{itemize}
\item[(i)] If $\mathcal{A}$ is an $R_0$-tensor, then Condition \ref{cond-1} holds.
\item[(ii)] If the solution set of TCP$(q, \mathcal{A})$ is nonempty and Condition \ref{cond-1} holds, then the solution set of TCP$(q, \mathcal{A})$ is bounded.
\item[(iii)] If the solution set of TCP$(q, \mathcal{A})$ is nonempty, then it is closed.
\end{itemize}
\end{Proposition}

{\bf Proof.} The proof of this proposition can be found in the proof of \cite[Theorem 3.1]{bhw-15}. We omit it here.  \ep

It follows from \cite[Theorem 3.2]{sq-15} that if $\mathcal{A}\in \mathbb{T}_{m,n}$ is an $R$-tensor, then the solution set of TCP$(q, \mathcal{A})$ is nonempty. Combining this with the fact that every $R$-tensor is an $R_0$-tensor, we can obtain the following result, which is a generalization of Theorem 3.2 in \cite{sq-15}.
\begin{Theorem}\label{thm-solu-r}
Suppose that $\mathcal{A}\in \mathbb{T}_{m,n}$ is an $R$-tensor and $q\in \mathbb{R}^n$ is given, then the solution set of TCP($q,\mathcal{A}$) is nonempty and compact.
\end{Theorem}

Now, we discuss properties of the solution set of TCP$(q, \mathcal{A})$ with ${\cal A}$ being an $ER$-tensor.

\begin{Theorem}\label{thm-solu}
Suppose that $\mathcal{A}\in \mathbb{T}_{m,n}$ is an $ER$-tensor and $q\in \mathbb{R}^n$ is given. Then the solution set of TCP($q,\mathcal{A}$) is nonempty and compact.
\end{Theorem}

{\bf Proof.} First, we show that the solution set of TCP($q,\mathcal{A}$) is nonempty. Suppose that TCP($q,\mathcal{A}$) has no solution. Then it follows from Theorem \ref{thm-ex-basic} that there exists an exceptional family of elements for $f(x)={\cal A}x^{m-1}+q$, i.e., there exists a sequence $\{x^k\}\subset \mathbb{R}^n_+$ satisfying $\|x^k\|\rightarrow \infty$ as $k\rightarrow \infty$ and, for each $k>0$, there exists a scalar $\mu_k>0$ such that
\begin{eqnarray}
({\cal A}x^k)^{m-1}_i+q_i&=&-\mu_kx_i^k,\quad \mbox{\rm if}\; x_i^k>0,\label{E-ex-f-1}\\
({\cal A}x^k)^{m-1}_i+q_i&\geq&0,\qquad\quad\; \mbox{\rm if}\; x_i^k=0. \label{E-ex-f-2}
\end{eqnarray}
Without loss of generality, we assume that $\frac{x^k}{\|x^k\|}\rightarrow x^*$. Then we have
\begin{equation}\label{E-ex-f-3}
x^*\in \mathbb{R}^n_+\quad \mbox{\rm and}\quad x^*\neq 0.
\end{equation}
We consider the following two cases:
\begin{itemize}
  \item Suppose that $i\in \{j\in [n]: x^*_j>0\}$, then $x_i^k>0$ for all sufficiently large $k$. Denote $t^*:=\lim_{k\rightarrow \infty}\frac{\mu_k}{\|x^k\|^{m-2}}$, then we have $t^*\geq 0$ and
\begin{eqnarray*}
-\frac{[{\cal A}(x^*)^{m-1}]_i}{x^*_i}
&=& \lim_{k\rightarrow \infty}\left[-\left({\cal A}\left(\frac{x^k}{\|x^k\|}\right)^{m-1}\right)_i\frac{1}{x^k_i/\|x^k\|}\right]\\
&=& \lim_{k\rightarrow \infty} \frac{-\left[({\cal A}(x^k)^{m-1})_i+q_i\right]}{\|x^k\|^{m-1}} \frac{1}{x^k_i/\|x^k\|}\\
&=& \lim_{k\rightarrow \infty} \frac{\mu_k x_i^k}{\|x^k\|^{m-1}}\frac{1}{x^k_i/\|x^k\|}\qquad (\mbox{\rm by}\; (\ref{E-ex-f-1}))\\
&=& t^*,
\end{eqnarray*}
which yields that $[{\cal A}(x^*)^{m-1}]_i+t^*x^*_i=0$ for any $i\in \{j\in [n]: x^*_j>0\}$.
  \item Suppose that $i\in \{j\in [n]: x^*_j=0\}$, then $x_i^k/\|x^k\|\rightarrow 0$ as $k\rightarrow \infty$; and hence, $\mu_k x_i^k/\|x^k\|\rightarrow 0$ as $k\rightarrow \infty$. Furthermore, we have
\begin{eqnarray*}
[{\cal A}(x^*)^{m-1}]_i&=&\lim_{k\rightarrow \infty}\left[{\cal A}\left(\frac{x^k}{\|x^k\|}\right)^{m-1}\right]_i =\lim_{k\rightarrow \infty}\frac{[{\cal A}(x^k)^{m-1}]_i+q_i}{\|x^k\|^{m-1}}\\
&=& \left\{\begin{array}{ll}
\lim\limits_{k\rightarrow \infty}\frac{[{\cal A}(x^k)^{m-1}]_i+q_i}{\|x^k\|^{m-1}}\geq 0,\quad &\mbox{\rm if}\; x_i^k=0,
\quad (\mbox{\rm by}\; (\ref{E-ex-f-2}))\vspace{2mm}\\
\lim\limits_{k\rightarrow \infty} \frac{-\mu_k x_i^k}{\|x^k\|^{m-1}}=0,\qquad & \mbox{\rm if}\; x_i^k>0,
\quad (\mbox{\rm by}\; (\ref{E-ex-f-1}))
\end{array}\right.
\end{eqnarray*}
which yields that $[{\cal A}(x^*)^{m-1}]_i\geq 0$ for any $i\in \{j\in [n]: x^*_j=0\}$.
\end{itemize}
These, together with (\ref{E-ex-f-3}), imply that $(x^*,t^*)\in(\mathbb{R}^n_+\setminus \{0\})\times\mathbb{R}_+$ satisfies the system (\ref{ertensor}), which contradicts the condition that ${\cal A}$ is an $ER$-tensor. Thus, the solution set of TCP($q,\mathcal{A}$) is nonempty.

Second, we show that the solution set of TCP($q,\mathcal{A}$) is compact. On one hand, since ${\cal A}$ is an $ER$-tensor, by Proposition \ref{pro1}(i) it follows that ${\cal A}$ is an $R_0$-tensor. This, together with Proposition \ref{prop-basic-1}(i)(ii), implies that the solution set of TCP($q,\mathcal{A}$) is bounded. On the other hand, by Proposition \ref{prop-basic-1}(iii), we have that the solution set of TCP($q,\mathcal{A}$) is closed. Therefore, we obtain that the solution set of TCP($q,\mathcal{A}$) is compact.

The proof is complete.\ep

A function $g: \mathbb{R}^n\rightarrow \mathbb{R}^n$ is called to be positively homogeneous of degree $t$ with $t$ being a positive integer if $g(\lambda x)=\lambda^tg(x)$ for all $\lambda>0$; and it is called to be positively homogeneous when $t=1$. It has been proved in \cite[Theorem 4.1]{zi-00} that if $g(x)=f(x)-f(0)$ is positively homogeneous and exceptionally regular, then CP$(f)$ has a solution. It should be pointed out that the proof of the existence of solution to TCP($q,\mathcal{A}$) given in Theorem \ref{thm-solu} is similar to the one of  \cite[Theorem 4.1]{zi-00}, however, the function $g(x)={\cal A}x^{m-1}$ involved in Theorem \ref{thm-solu} is  positively homogeneous of degree $m-1$.
\begin{Remark}
Since a positive definite tensor is an $ER$-tensor and so is a strictly copositive tensor, the results in \cite[Theorem 4.5]{cqw} are special cases of Theorem \ref{thm-solu}. In addition, by Theorem \ref{wpiser-1}, the class of $wP$-tensors is a proper subset of the class of $ER$-tensors, and hence, the result of \cite[Theorem 6.2]{dlq-15} is also a special case of Theorem \ref{thm-solu}.
\end{Remark}

From Theorem \ref{thm-solu}, we know that for any $q\in \mathbb{R}^n$, TCP($q, \mathcal{A}$) with $\mathcal{A}$ being an $ER$-tensor has at least one solution. So, we have the following corollary.
\begin{Corollary}
Every $ER$-tensor is a $Q$-tensor.
\end{Corollary}

From Theorem \ref{thm-solu} and Corollary \ref{cor-thm-equi}, we have the following result.

\begin{Corollary}\label{cor-thm-solu}
If $\mathcal{A}\in \mathbb{T}_{m,n}$ is a $P_0+R_0$-tensor, then the solution set of TCP($q,\mathcal{A}$) is nonempty and compact.
\end{Corollary}

\section{Conclusions}
In this paper, we introduced the concept of $ER$-tensor. We showed that many important classes of tensors are the subclasses of the class of $ER$-tensors, including the class of $wP$-tensors which is recently defined by Ding, Luo and Qi\cite{dlq-15}. We also proved that the intersection of the class of $ER$-tensors and the class of $R$-tensors is nonempty, and by constructing two examples, we showed that an $R$-tensor is not always an $ER$-tensor and an $ER$-tensor is not always an $R$-tensor. We studied some properties of $ER$-tensor, in particular, we obtained that the equivalence of the class of $R_0$-tensors, the class of $R$-tensors and the class of $ER$-tensors within the class of the semi-positive tensors. As an application, we investigated the tensor complementarity problem with the involved tensor being an $ER$-tensor and showed that the solution set of this class of tensor complementarity problems is nonempty and compact. We also obtained that the solution set of the tensor complementarity problem, with the involved tensor being an $R$-tensor or a $P_0+R_0$-tensor, is nonempty and compact. We believe that more properties related to $ER$-tensor can be further studied in the future work. In addition, it is worth investigating the properties and applications of $ER$-matrix.

\end{document}